\documentclass[colorlinks]{goose-article}

\usepackage{multirow, booktabs}
\usepackage{subcaption}
\usepackage[capitalise]{cleveref}

\author[1]{Jasper~C.\ Volmer}
\author[2]{Tom~W.J.\ de Geus}
\author[1]{Ron~H.J.\ Peerlings}

\affil[1]{%
  Department of Mechanical Engineering,
  Eindhoven University of Technology,\nl
  P.O.\ Box 513,
  5600 MB Eindhoven,
  The Netherlands}

\affil[2]{%
  Institute of Physics,
  \'{E}cole Polytechnique F\'{e}d\'{e}rale de Lausanne (EPFL),\nl
  CH-1015 Lausanne,
  Switzerland}

\date{}

\title{Improving the initial guess for the Newton-Raphson protocol in time-dependent simulations}

\header{Journal of Computational Physics, \doi{10.1016/j.jcp.2020.109721}, \arxivid{1912.12140}}

\begin{document}

\maketitle

\keywords{
    Visco-plasticity;
    Newton-Raphson;
    Linearisation;
    FFT solver;
    Backward Euler}

\section{Introduction}

In simulating real-world problems using quasi-static mechanics,
one often models the material's constitutive response using a strain-rate dependent law.
One naturally does this when the material possesses a time dependent relaxation mechanism.
But it is also common to introduce an artificial strain-rate dependency to
regularise numerical solutions of a rate-independent material,
e.g.\ in (crystal) plasticity or damage simulations.

In this setting, we most frequently solve the balance of linear momentum in the following form
\begin{equation}
    \label{eq:equilibrium}
    \nabla \cdot
    \bm{\sigma} \left( \bm{\varepsilon}, \dot{\bm{\varepsilon}}, t, \ldots \right) =
    \bm{0}
\end{equation}
i.e.\ the divergence of the stress $\bm{\sigma}$ has to vanish everywhere in the domain
(see \cref{app:nomenclature} for our nomenclature).
This problem is generally hard to solve because of the complexity of the stress response
at the material point level,
as it often depends non-linearly on the strain $\bm{\varepsilon}$,
the strain rate $\dot{\bm{\varepsilon}}$,
the time $t$, and the loading history (carried in a number of internal variables).
This partial differential equation thus has to be solved numerically.
To this end, one typically discretises space,
resulting in a system of (non-)linear algebraic equations
\begin{equation}
    \label{eq:nonlinear-system}
    \underline{\bm{f}}
    \left\{
        \underline{\bm{\sigma}} ( \underline{\bm{\varepsilon}},
        \dot{\underline{\bm{\varepsilon}}},
        t,
        \ldots )
    \right\} = \underline{\bm{0}}
\end{equation}
where $\underline{f} \{ \underline{\bullet} \}$ is a linear or non-linear functional related
to the discretisation in space, whereby $\underline{\bullet}$ refers to
a discrete set of variables.

By and large, the most popular protocol for solving the system of non-linear equations in
Eq.~\eqref{eq:nonlinear-system} is the Newton-Raphson procedure.
This procedure employs a first-order approximation of Eq.~\eqref{eq:nonlinear-system} in
the neighbourhood of an approximate solution $\underline{\bm{\varepsilon}}_{i}$,
resulting in a linear system of the form
\begin{equation}
    \label{eq:lin_approx}
    \left.
        \frac{\partial \underline{\bm{f}}}{\partial \underline{\bm{\varepsilon}}}
    \right|_{\underline{\bm{\varepsilon}}_{i}}
    \delta \underline{\bm{\varepsilon}} =
    - \underline{\bm{f}}(\underline{\bm{\varepsilon}}_{i}, \ldots)
\end{equation}
that can be solved for $\delta \underline{\bm{\varepsilon}}$.
The approximate solution is then updated according to
\begin{equation}
    \label{eq:epsilon_update}
    \underline{\bm{\varepsilon}}_{i+1} =
      \underline{\bm{\varepsilon}}_{i} +
      \delta \underline{\bm{\varepsilon}}
\end{equation}
These so-called iterations are repeated until Eq.~\eqref{eq:nonlinear-system}
is satisfied with sufficient precision.
In many cases $\underline{\bm{\sigma}}$ contains ordinary differential equations
in time to describe the evolution of the internal variables.
These hidden ordinary differential equations are solved by discretising time,
often by some finite difference scheme.
As a result, Eq.~\eqref{eq:lin_approx} is employed consecutively at different points in time.

The computational efficiency of such a scheme relies crucially on
i) the accuracy of the first order approximation in Eq.~\eqref{eq:lin_approx} \cite{Simo_1985} and
ii) the accuracy of the initial guess $\bm{\varepsilon}_{0}$ that is iteratively refined using Eq.~\eqref{eq:epsilon_update}.
i) requires a consistent tangent at the material point level and its derivation is usually
well established \cite{Ju_1990, Souza_2011, Simo_2006}.
For ii) the obvious choice is to take the last known converged state as the starting point
$\bm{\varepsilon}_{0}$, however,
we show that for time-dependent problems a better choice can be made.
It involves a subtle interaction between the non-linear solver and the time dependence.
This interaction becomes obvious in the derivation of
i) by properly linearising all terms that are part of the discrete time integration scheme.
More specifically, we show that a step along the discrete time axis will lead to a viscous flow,
regardless of how the system is driven.
This is incorporated by an additional stress (or force) term present only in the first iteration
after the time increment amending to a logical choice for ii).

The purpose of this note is to present a derivation that naturally leads to the additional
term for the first iteration of a new time increment.
Furthermore, we show that the additional term can be easily interpreted as an initial guess
for the Newton-Raphson protocol.
We benchmark the improvement by solving the mechanical response of a dual-phase steel
microstructure using a modern numerical method based on the Fast Fourier Transform (FFT)
\cite{Zeman_2017, GooseFFT, Geus_2017}.
A reduction of the computation time of around 45\% is observed in comparison to taking the
last known converged state as an initial guess.
We emphasise that we present the procedure on a relatively simple model,
but that it applies to more complex models as well.
We illustrate this by considering also a more involved time integration scheme,
and thereby show that the procedure is straightforward to apply.

The remainder of this note is structured as follows:
A visco-plastic (time-dependent) material model is introduced together with its linearisation.
We thereby distinguish two components: the classical consistent tangent used in every iteration,
and an additional driving force inserted during the first iteration.
The performance of the classical and improved schemes is examined lastly.

\section{Material model}

\begin{figure}[htp]
    \centering
    \begin{subfigure}[t]{0.45\linewidth}
        \centering
        \includegraphics[width=1.\linewidth]{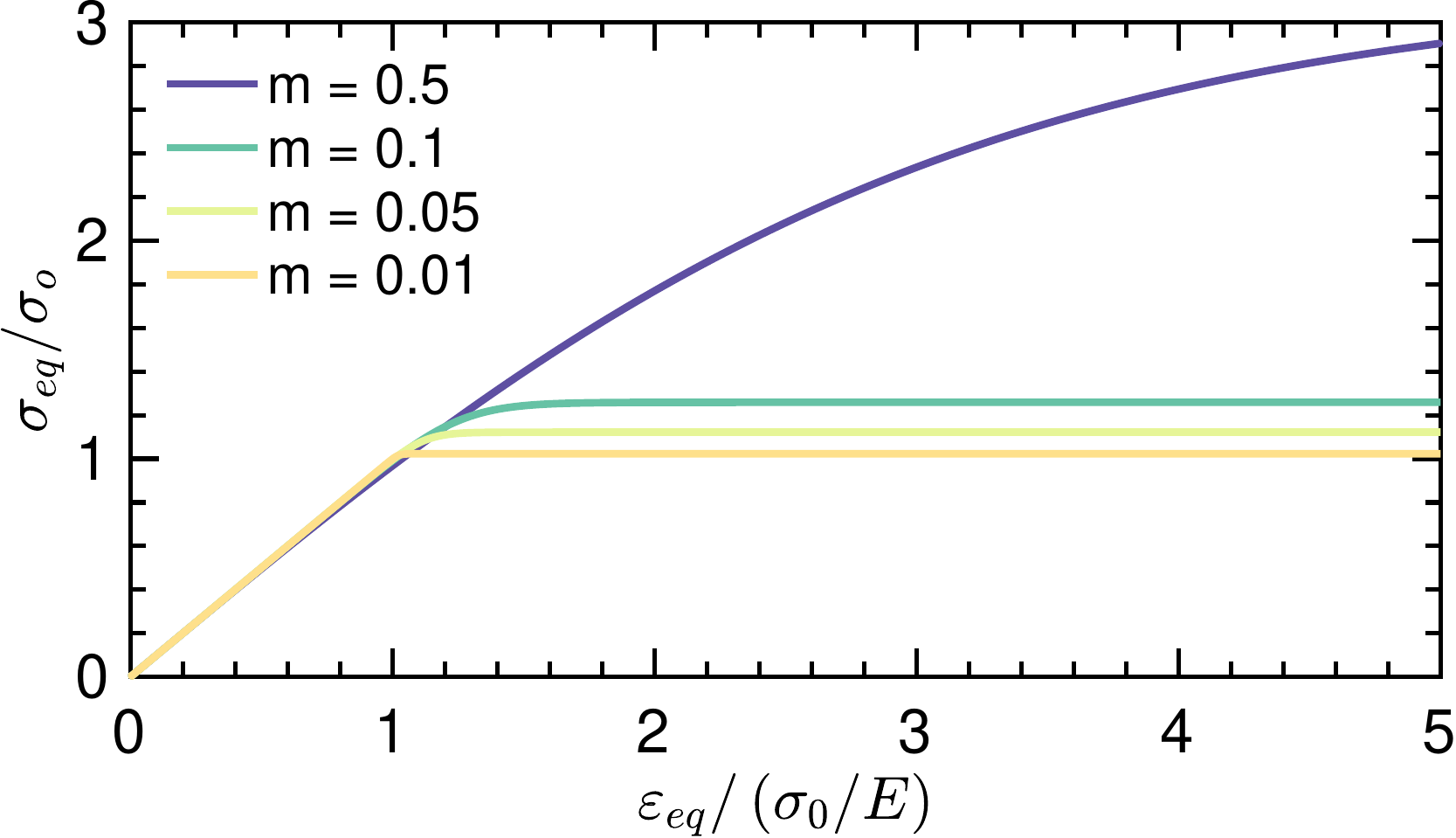}
        \caption{Influence of the rate sensitivity exponent $m$.}
        \label{fig:behaviour_m}
    \end{subfigure}\hspace{.05\linewidth}
    \begin{subfigure}[t]{0.45\linewidth}
        \centering
        \includegraphics[width=1.\linewidth]{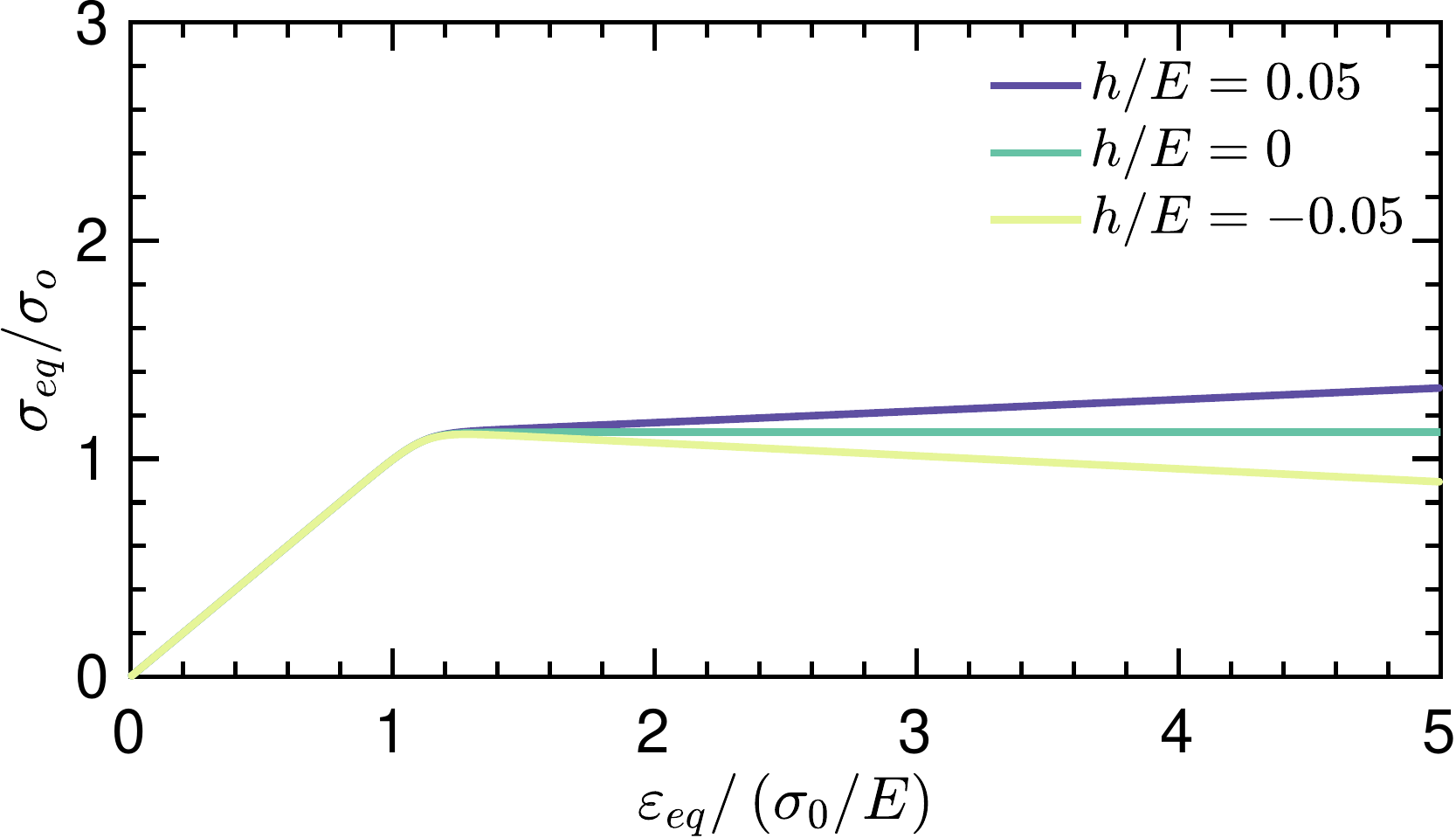}
        \caption{Influence of the hardening modulus $h$.}
        \label{fig:behaviour_h}
    \end{subfigure}
    \caption{
        Normalised stress-strain curves for an individual visco-plastic material point
        for $\sigma_0 / E = 0.01$ and an applied strain rate of
        $\dot{\varepsilon}_{eq}/{\dot {\gamma}_0 } = 10$ in the normal direction.
        The equivalent strain rate is defined work conjugate to the equivalent stress:
        ${\varepsilon_{eq}} \equiv \sqrt{2/3 \, \bm{\varepsilon}^d:\bm{\varepsilon}^d}$.}
    \label{fig:behaviour_vp}
\end{figure}

A relatively simple visco-plastic model, based on the small strain assumption, is used here.
The employed model is described in \cite[chapter 11]{Souza_2011}.
The stress is set by the elastic strain using Hooke's law.
Thereto, the total strain $\bm{\varepsilon}$ is additively split into an elastic part
$\bm{\varepsilon}_{e}$ and a plastic part $\bm{\varepsilon}_{p}$ as
\begin{equation}
    \label{eq:strain_split_3D}
    \bm{\varepsilon} \equiv \bm{\varepsilon}_{e} + \bm{\varepsilon}_{p}
\end{equation}
and thus
\begin{equation}
    \label{eq:sigma_3D}
    \bm{\sigma} \equiv {}^4\bm{C}_e:\,{\bm{\varepsilon}_e}
\end{equation}
with
\begin{equation}
    \label{eq:4ce}
    {}^4\bm{C}_e \equiv K\,\bm{I} \bm{I} + 2G\,\,{{}^4 \bm{I^d}}
\end{equation}
where the fourth-order tensor ${}^4\bm{I}^d \equiv {}^4\bm{I}^s - \bm{I} \bm{I} /3$
projects an arbitrary tensor $\bm{A}$ onto its deviatoric part
$\bm{A}^{d} = \bm{A} - \text{tr}(\bm{A})\bm{I}/3$.
The elastic material parameters are the bulk modulus $K$ and the shear modulus $G$,
which depend on Young's modulus $E$ and Poisson's ratio $\nu$ in the usual way.
The evolution of plastic strain is given by the flow rule as
\begin{equation}
    \label{eq:epspdot}
    \dot{\bm{\varepsilon}}_{p}
    \equiv \dot{\gamma} \, \bm{N}
\end{equation}
whereby the direction of the plastic flow is given by
\begin{equation}
    \label{eq:N}
    \bm{N} \equiv \frac{3}{2}\frac{{{\bm{\sigma}^d}}}{{{\sigma_{eq}}}}
\end{equation}
where $\bm{\sigma}^d$ is the deviatoric part of the stress
($\bm{\sigma}^d = {}^4 \bm{I}^d : \bm{\sigma}$) and ${\sigma_{eq}}$
is the Von Mises equivalent stress
(${\sigma_{eq}} \equiv \sqrt{3/2 \, \bm{\sigma}^d:\bm{\sigma}^d}$).
The magnitude of the plastic flow $\dot{\gamma}$ is given by  Norton's rule as

\begin{equation}
    \label{eq:gdot}
    \dot \gamma \equiv
    {\dot \gamma_0}{\left( {\frac{{{\sigma_{eq}}}}{{{\sigma_s}}}} \right)^{1/m}}
\end{equation}
The material constants are:
the reference strain rate ${\dot \gamma_0}$,
the rate-sensitivity exponent $m$ and
the yield stress of the material $\sigma_s$.
Note that $\dot{\gamma}$ is by construction nonnegative,
and depends non-linearly on the stress $\bm{\sigma}$,
and therefore on the plastic strain $\bm{\varepsilon}_p$,
through the rate sensitivity exponent $m$.
We, furthermore, let the yield stress $\sigma_s$ evolve with
the accumulated plastic strain as follows:
\begin{equation}
    \label{eq:hardening_law}
    {\sigma_s} \equiv {\sigma_o} + h\,{\varepsilon_p}
\end{equation}
where $\sigma_0$ is the initial yield stress and $h$ is the hardening modulus.
If $h=0$ then the model behaves perfectly plastically,
whereas it hardens when $h>0$ and softens when $h < 0$.
Finally, the accumulated plastic strain ${\varepsilon_p}$ is determined from
\begin{equation}
    \label{eq:plastic-strain}
    {\varepsilon_p} \equiv \int\limits_0^t {\dot \gamma dt'}
\end{equation}
To illustrate the behaviour of the visco-plastic model introduced above,
several normalised stress-strain curves for a single material point are presented in
Fig.~\ref{fig:behaviour_vp}.
Fig.~\ref{fig:behaviour_m} shows the behaviour of the visco-plastic model for several
values of the rate sensitivity exponent $m$.
A sharp transition from the elastic to the plastic regime,
as would be observed for rate-independent elasto-plastic behaviour,
can be approximated by a small value for the rate sensitivity exponent $m$.
The different regimes of hardening,
perfect plasticity and softening are shown in Fig.~\ref{fig:behaviour_h}.

\section{Time discretisation}
\label{sec:time_discretisation}

The numerical treatment proceeds by discretising the material model in time,
for which several choices exist.
We employ the implicit Backward Euler protocol,
known for its ease of implementation and robust convergence.
The procedure is however general, which we demonstrate using the
more advanced generalised trapezoidal integration scheme in \ref{app:theta_integration}.

Using the Backward Euler protocol,
the discrete version of the flow rule in Eq.~\eqref{eq:epspdot} reads
\begin{equation}
    \label{eq:Delta_epsp}
    \Delta \bm{\varepsilon}_{p} = \Delta \gamma \; \bm{N}^{t+\Delta t}
\end{equation}
where
\begin{equation}
    \label{eq:BE_dgam}
    \Delta \gamma =
    \gamma^{t+\Delta t} - \gamma^{t} = \Delta t \, {{\dot \gamma }^{t + \Delta t}}
\end{equation}
To compute the unknowns $\gamma^{t+\Delta t}$ and $\bm{N}^{t+\Delta t}$,
we apply the common approach of performing a radial return map.
This involves formulating a trial state ${}^{tr}{\bullet}$,
in which a strain increment is assumed fully elastic,
from which the plastic evolution is sought to end up in an admissible state.
For a model like the one presented here, this involves solving a single scalar,
yet non-linear, equation, as it is easily shown that the stress directions in
the trial state are the same as that of the admissible state.
We may therefore write $\bm{N}^{t+\Delta t} = {}^{tr}\bm{N}$,
see e.g.\ \cite{Souza_2011} for details.

\section{Consistent linearisation of the stress update around an arbitrary reference state}

The consistent tangent relates, to the first order,
a perturbation in the strain, $\delta \bm{\varepsilon}$,
to the resulting perturbation in the stress, $\delta \bm{\sigma}$.
In this section we derive this tangent from multivariable linearisation around a
fully known `reference state'.
The reference state is left undefined for the moment, but will be clarified below.
This approach allows us to clearly identify the terms that are proposed as improvement.

The procedure begins by writing all unknown variables at time $t+\Delta t$ as a fully known
reference state~${\bullet^*}$ plus a small perturbation $\delta \bullet$.
In terms of our model we have:

\begin{equation}
    \label{eq:sigma_NR_3D}
    {\bm{\sigma}^{t + \Delta t}} \equiv {\bm{\sigma}^*} + \delta \bm{\sigma}
\end{equation}
\begin{equation}
    \label{eq:epsilon_NR_3D}
    {\bm{\varepsilon}^{t + \Delta t}} \equiv {\bm{\varepsilon}^*} + \delta \bm{\varepsilon}
\end{equation}
\begin{equation}
    \label{eq:epsilon_p_NR_3D}
    \bm{\varepsilon}_p^{t + \Delta t} \equiv \bm{\varepsilon}_p^* + \delta {\bm{\varepsilon}_p}
\end{equation}
\begin{equation}
    \label{eq:N_lin}
    {\bm{N}^{t + \Delta t}} \equiv {\bm{N}^*} + \delta \bm{N}
\end{equation}
\begin{equation}
    \label{eq:gamma_lin}
    {\gamma^{t + \Delta t}} \equiv {\gamma^*} + \delta \gamma
\end{equation}
\begin{equation}
    \label{eq:gamma_dot_lin}
    {\dot{\gamma}^{t + \Delta t}} \equiv {\dot{\gamma }^*} + \delta \dot{\gamma }
\end{equation}
The next step is to linearise all the equations used in the calculation of
the constitutive response around the reference state ${\bullet^*}$.
The elastic law in Eq.~\eqref{eq:sigma_3D} is already linear and hence gives
\begin{equation}
    \label{eq:lin_constit_RI_3D}
    \delta \bm{\sigma} =
    {}^4\bm{C}_e:\left( {\delta \bm{\varepsilon} - \delta {\bm{\varepsilon}_p}} \right)
\end{equation}
A first-order approximation of $\delta {\bm{\varepsilon}_p}$ is slightly more involved
as it is non-linear and time-dependent.
Its derivation starts from Eq.~\eqref{eq:Delta_epsp} and eventually leads to

\begin{equation}
    \label{eq:del_epsp}
    \delta {\bm{\varepsilon}_p}
    = \delta \gamma \,{\bm{N}^*}
    + \underbrace{\left({\gamma^*}
    - {\gamma^{* - \Delta t}}\right)}_{\Delta {\gamma^*}} \delta {\bm{N}}
\end{equation}
Note that $\Delta \gamma^*$ and $\bm{N}^*$ are known quantities.
The derivation continues by developing expressions for $\delta \gamma$ and $\delta {\bm{N}}$,
to acquire a closed-form expression for the small variation of plastic strain
$\delta {\bm{\varepsilon}_p}$.
To obtain $\delta{\bm{N}}$ we use the result from the radial return map,
so that $\delta{\bm{N}} = \delta({}^{tr}{\bm{N}})$.
The latter can be entirely evaluated in the trial state, and results in
\begin{equation}
    \label{eq:del_N}
    \delta \bm{N} =
    \left[ {\frac{{3G}}{{{}^{tr}\sigma_{eq}^*}}{}^4{\bm{I}^d} -
    \frac{{2G}}{{{}^{tr}\sigma_{eq}^*}}{\bm{N}^*}{\bm{N}^*}} \right] :
    \delta \bm{\varepsilon}
\end{equation}
To find $\delta \gamma$, we combine the results of Eqs.~\eqref{eq:BE_dgam},
\eqref{eq:gamma_lin} and \eqref{eq:gamma_dot_lin} into
\begin{equation}
    \label{eq:lin_BE}
    \gamma^* - \gamma^t - \Delta t {\dot \gamma}^{*} + \delta \gamma
    - \Delta t \delta\dot \gamma = 0
\end{equation}
By linearising Eq.~\eqref{eq:gdot} around ${\bullet^*}$,
the small variation $\delta\dot \gamma$ can be written as
\begin{equation}
    \label{eq:lin_del_gam}
    {\delta\dot \gamma } =
    \frac{{\partial \dot \gamma }}{{\partial {\sigma_{eq}}}}\delta {\sigma_{eq}} +
    \frac{{\partial \dot \gamma }}{{\partial {\sigma_s}}}\delta {\sigma_s} =
    \frac{\alpha^*}{\Delta t} \left( {2{\bm{N}^*}:\delta \bm{\varepsilon } - 3\delta \gamma } -
    \frac{{\sigma_{eq}^*\,h\,\delta \gamma }}{{\sigma_s^*}G} \right),
    \qquad
    {\alpha^*} =
    \frac{{{{\dot \gamma }_0 G \Delta t}}}{{m\sigma_s^* }}{
    \left( {\frac{{\sigma_{eq}^*}}{{\sigma_s^*}}} \right)^{\frac{1}{m} - 1}}
\end{equation}
Thereby we have employed the results of Eqs.~\eqref{eq:lin_constit_RI_3D},
\eqref{eq:del_epsp} and \eqref{eq:del_N}.

A closed-form expression for $\delta \gamma$ can now be established by substitution of
Eq.~\eqref{eq:lin_del_gam} in Eq.~\eqref{eq:lin_BE}.
Substituting that expression for $\delta \gamma$ in Eq.~\eqref{eq:del_epsp} and
the resulting expression in Eq.~\eqref{eq:lin_constit_RI_3D}
finally gives the consistent linearisation of the stress update, as follows
\begin{equation}
    \label{eq:del_sigma_a}
    \delta \bm{\sigma} =
    {}^4{\bm{C}_{vp}^*} : {\delta \bm{\varepsilon} -
    \frac{G \beta^*}{\alpha^*}{\left(
    {{\Delta t{{\dot \gamma }^*} - \gamma^* + \gamma^t}} \right)}{\bm{N}^*}}
\end{equation}
where the consistent tangent ${}^4{\bm{C}_{vp}^*}$ for the visco-plastic model reads
\begin{equation}
    \label{eq:gen_C_vp}
    {}^4{\bm{C}_{vp}^*} =
    {}^4{\bm{C}_e} -
    2G\beta^*{\bm{N}^*}{\bm{N}^*} -
    \frac{\Delta {\gamma^*} 4G^2}{^{tr}\sigma_{eq}^*}
    \left[ { \frac{3}{2}{}^4{\bm{I}^d} - {\bm{N}^*}{\bm{N}^*}} \right]
\end{equation}
and
\begin{equation}
    \label{eq:beta}
    \beta^* = \frac{{2{\alpha^*}}}{{1 + 3{\alpha^*} +
    \frac{{\sigma_{eq}^*}h}{{\sigma_s^*}G} \alpha^*}}
\end{equation}
which can be further reorganised\footnote{
    It may be helpful to realise that the following identity holds:
    $\frac{G \beta^*}{\alpha^*} = \kappa^* \left( 2G - 3G\beta^* \right)$}
to:
\begin{equation}
    \label{eq:del_sigma_b}
    \delta \bm{\sigma} =
    {}^4{\bm{C}_{vp}^*} : \left[ {\delta \bm{\varepsilon} -
    {\left( {{\Delta t{{\dot \gamma }^*} - \gamma^* + \gamma^t}} \right)}
    \kappa^*{\bm{N}^*}} \right]
\end{equation}
with
\begin{equation}
    \label{eq:kappa}
    \kappa^* =
    \frac{1}{\left( 1 + \frac{{\sigma_{eq}^*}h}{{\sigma_{s}^*} G} \alpha^* \right)}
\end{equation}

Note how the choice of the reference state ${\bullet^*}$ determines at which state
the consistent tangent is evaluated and that it does not affect the expression itself.
It does, however, affect the relevance of the second term between brackets in
Eq.~\eqref{eq:del_sigma_b}, as we will see next.

\section{Reference state}
%
\subsection{Recovering the classic Newton-Raphson iteration}

We now define the reference state denoted by ${\bullet^*}$.
We first consider `ordinary'
Newton-Raphson iterations within one discrete time step ${\bullet^{t + \Delta t}}$
as for example in Eq.~\eqref{eq:epsilon_update}.
In this case, an iterative update of the unknown(s) is obtained by linearising around the
last known iterative state, denoted by the iteration counter $i$.
In this case, Eq.~\eqref{eq:del_sigma_b} reduces to the classical
\cite{Ju_1990, Souza_2011, Simo_2006}:
\begin{equation}
    \label{eq:del_sigma_NR}
    \delta \bm{\sigma} = {}^4{\bm{C}_{vp}^{i}}:{\delta \bm{\varepsilon}}
\end{equation}
where the reference state ${\bullet^*} \equiv {\bullet^{i}}$ at $t + \Delta t$,
the latter not being explicitly included in the notation.
This result follows from Eq.~\eqref{eq:del_sigma_b} as
$\gamma^* \equiv \gamma^i$ and $\dot{\gamma}^* \equiv \dot{\gamma}^i$.
Recognising the discretised strain rate ($\Delta t \dot{\gamma}^i \equiv \gamma^i - \gamma^t $)
we thus find the three rightmost terms in Eq.~\eqref{eq:del_sigma_b} to cancel.

\subsection{Obtaining the improved initial guess}

For the first iteration of every new time increment, we have to be careful.
Commonly, one simply uses the last available tangent as in Eq.~\eqref{eq:del_sigma_NR}.
This would amend to taking the tangent of the last iteration $i$ of the previous time step
(at time $t$), that resulted in a converged state.
We argue that when taking this converged state as our reference state for linearisation
(${\bullet^*} \equiv {\bullet^{t}}$), an extra term appears in the stress update:
\begin{equation}
    \label{eq:del_sigma_load}
    \delta \bm{\sigma} =
    {}^4{\bm{C}_{vp}^{t}}: \left[ {\delta \bm{\varepsilon}} -
    {\Delta t \, {{\dot \gamma }^t}}\kappa^t{\bm{N}^t} \right]
\end{equation}
where ${}^4{\bm{C}_{vp}^{t}}$ is the consistent tangent according to
Eq.~\eqref{eq:gen_C_vp} evaluated at the converged state at time $t$.
Note that this result trivially follows from Eq.~\eqref{eq:del_sigma_b} as
$\gamma^* \equiv \gamma^t$ and the two rightmost terms cancel.
The extra term in Eq.~\eqref{eq:del_sigma_load} (cf.\ Eq.~\eqref{eq:del_sigma_NR})
can be interpreted as the increase in plastic strain $\bm{\varepsilon}_p$
over the time step $\Delta t$ as caused by the stress $\bm{\sigma}$ at time $t$.
Note that the magnitude of this plastic strain increase computed from
$ \Delta t \, {{\dot \gamma }^t}$ is scaled with the variable $\kappa^t$,
which takes into account the effect of the plastic strain increase on the yield stress,
i.e.\ $\kappa^t > 1$ for hardening,
$\kappa^t  = 1$ for perfect plasticity and $\kappa^t < 1$ for softening.
Naturally, this expected increase in plastic strain, based on variables at time $t$,
is only an estimate.
As a result, $\delta \bm{\sigma}$ is only a prediction of the incremental change in stress.
This prediction thereby effectively sets an initial guess
from which to start the regular Newton-Raphson iterations.

The avid reader may wonder if for the first iteration after a time increment,
i.e.\ the situation described above, it would not simply suffice to use an explicit
increment to yield the same result.
It is emphasised that interchanging Eq.~\eqref{eq:BE_dgam} with an explicit substitute
(e.g.\ Forward Euler) yields:
i) a different expression for the consistent tangent ${}^4{\bm{C}_{vp}^{t}}$
in Eq.~\eqref{eq:del_sigma_load} and
ii) no compensation for the change in yield stress as the variable $\kappa^t$ does not appear.

To conclude, it is emphasised that the extra term in Eq.~\eqref{eq:del_sigma_load}
is a result of the load increment and the time-dependent material model.
It is therefore only included in the first iteration after a time increment.
In the regime where the plastic flow is negligible ($\dot {\underline{\gamma }}^t \approx 0$)
or for rate-independent material models, there is no contribution of the improved initial guess.

\section{Case study}

To show its relevance, we employ the improved initial guess in a case study.
Thereby we make use of a modern numerical solution procedure for micro-mechanical problems,
that is based on the Fast Fourier Transform (FFT).
As extensively described in \cite{Zeman_2017}, like in the Finite Element Method,
Eq.~\eqref{eq:equilibrium} is solved in a weak sense.
The resulting volume integral is evaluated numerically by introducing nodal unknowns that
are distributed on a regular grid (i.e.\ pixels or voxels).
Owing to this choice, they can be interpolated using globally supported trigonometric polynomials.
Numerical quadrature then proceeds by evaluating equally weighted quadrature points that
coincide with the nodes.
The result is a scheme in which essentially local equilibrium equations are coupled by
the application of the Fourier transform and its inverse, which can be done using efficient
and mature FFT libraries.
The details of how the improved initial guess appears in the algorithm proposed by
\cite{Zeman_2017} are given in \ref{app:implementation}.

\begin{figure}[tbp]
    \begin{subfigure}[t]{0.33\linewidth}
        \centering
        \includegraphics[width=.58\linewidth]{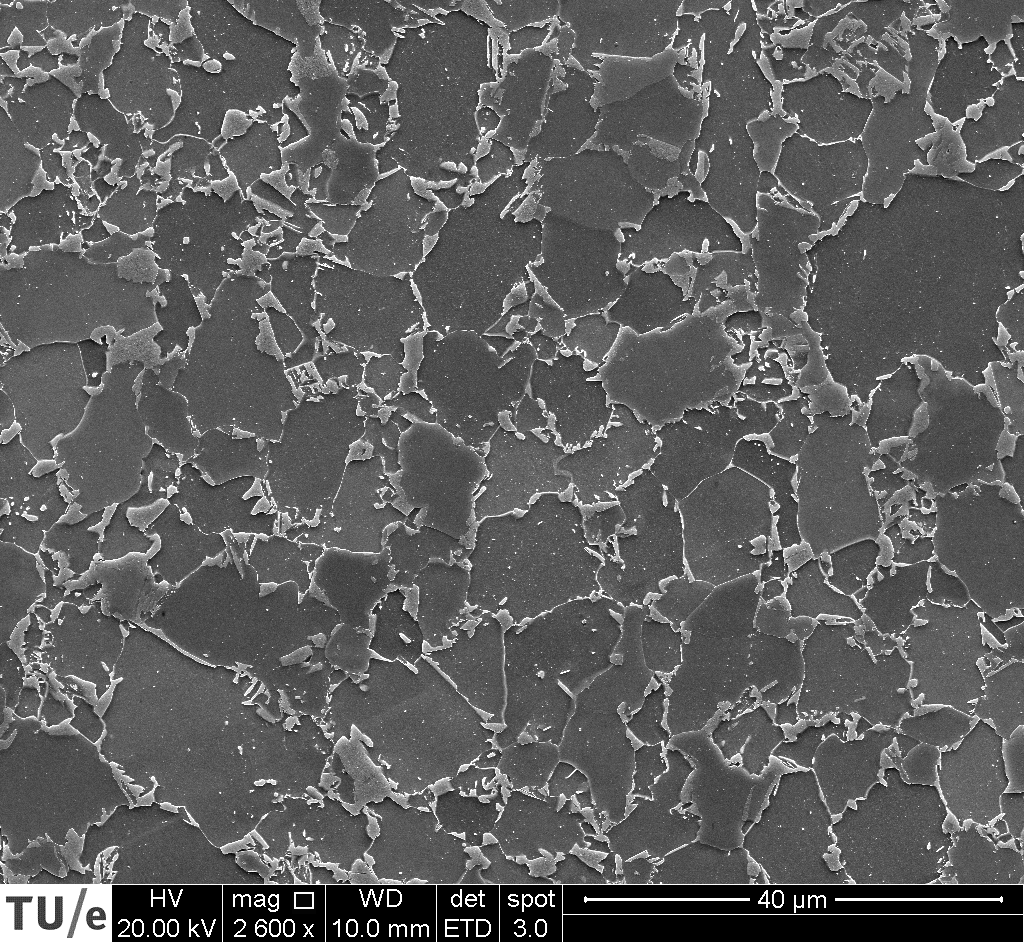}
        \caption{SEM image after scanning.}
        \label{fig:micrograph_a}
    \end{subfigure}
    \begin{subfigure}[t]{0.33\linewidth}
        \centering
        \includegraphics[width=.86\linewidth]{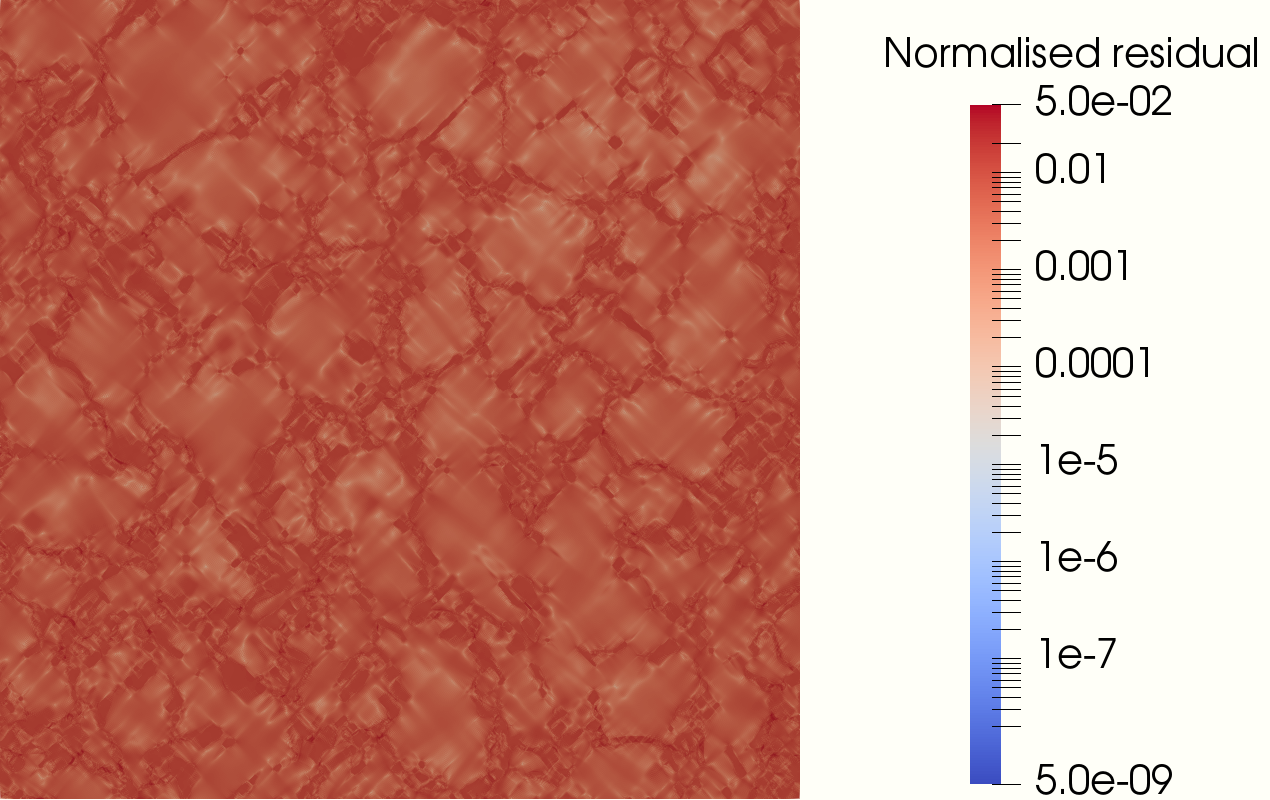}
        \caption{Initial guess as suggested in \cite{Zeman_2017}.}
        \label{fig:oig}
    \end{subfigure}
    \begin{subfigure}[t]{0.33\linewidth}
        \centering
        \includegraphics[width=.86\linewidth]{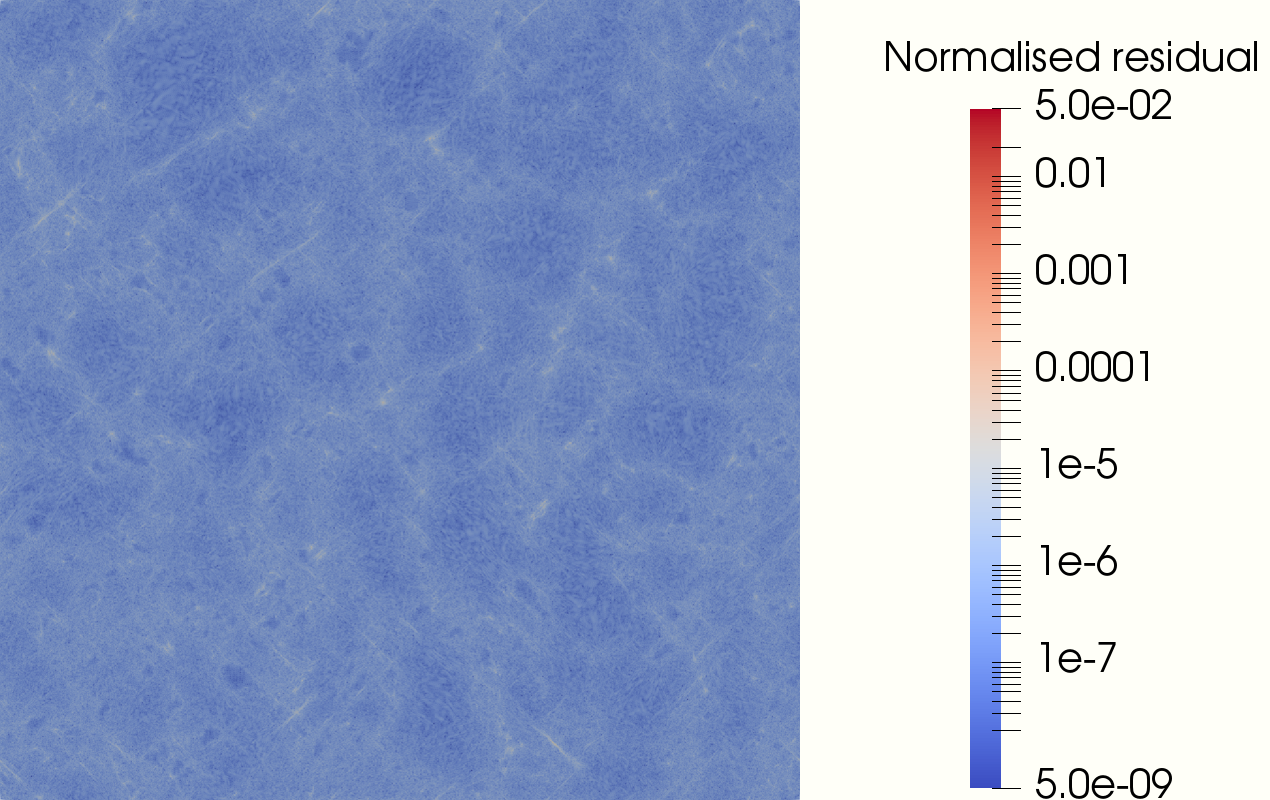}
        \caption{Initial guess using the extra term of Eq.~\eqref{eq:del_sigma_load}.}
        \label{fig:iig}
    \end{subfigure}
    \caption{
        SEM image (801 x 801 pixels) and visualisation of the magnitude of the residual of
        the mechanical equilibrium equation normalised with the yield stress of martensite at
        the start of the Newton-Raphson iterations on the SEM micrograph for
        perfect visco-plasticity at $\varepsilon_{appl} = 5.0 \cdot 10^{-2}$.}
    \label{fig:ig_comparisson}
\end{figure}

We study the efficiency of our improvement based on a realistic example in which we compute
the microscopic response of a microstructure that is subjected to a macroscopic shear strain.
The microstructure is taken from a micrograph of a commercial dual-phase steel sample (DP600),
acquired using a scanning electron microscope (SEM), as shown in Fig.~\ref{fig:micrograph_a}.
Dual-phase steel consists of two main constituents: i) ferrite, a soft and ductile phase,
which shows up in dark in the micrograph in Fig.~\ref{fig:micrograph_a} and ii) martensite,
a hard and brittle phase, which shows up bright in Fig.~\ref{fig:micrograph_a} and has
a volume fraction of about 17\% in this image of 801x801 pixels.
For our case study we assume that the microstructure is continuous and consists only of these
two phases, which we both assume to obey the visco-plastic model presented above.
To this end the micrograph in Fig.~\ref{fig:micrograph_a} is thresholded\footnote{
    Both the micrograph and the corresponding binary image obtained by thresholding have
    been taken from the GooseFFT repository \cite{GooseFFT, Geus_2017},
    see \cite{Geus_2016} for the experimental and thresholding protocol.}
to obtain a binary image.
Each pixel then corresponds to a nodal point for the FFT-solver,
whereby the material parameters are different depending on the phase,
see Table \ref{tbl:mat_par}.
Note that we consider three cases: hardening, perfect plasticity, and softening.
The parameters for these cases are loosely based on \cite{Sun_2009} and \cite{Eisenlohr_2013}.

In the simulations, the specimen is subject to periodic boundary conditions
(as so required by our solver, but common in this type of homogenisation problems).
An average strain $\bar{ \bm{\varepsilon}}$ is prescribed which induces a
pure shear strain according to
\begin{equation}
    \label{eq:eps_bc}
    \bar{ \bm{\varepsilon}}  =
    \frac{{\sqrt 3 }}{2}{{\varepsilon }_{appl}}
    \left( {{{\vec e}_x}{{\vec e}_x} - {{\vec e}_y}{{\vec e}_y}} \right)
\end{equation}
where $\varepsilon_{appl}$ is the applied strain and $\vec e_x$ and $\vec e_y$ are
the unit vectors, respectively in the horizontal and vertical direction.
For the simulations with hardening and perfect plasticity, the applied strain was
incrementally increased to $\varepsilon_{appl} = 0.05$ at a strain rate of
$\dot{\varepsilon}_{appl} = 0.01$ [1/s] in 100 time steps.
For the test cases that include softening, the equivalent strain was incrementally increased to
$\varepsilon_{appl} = 0.01$ with the same strain rate and number of time steps.

Figs.~\ref{fig:ig_comparisson}(b) and \ref{fig:ig_comparisson}(c) give an
example of the effect of
the extra term using the perfectly plastic material model.
The residual is visualised based on the computation of the mechanical equilibrium and
normalised with the yield stress of martensite.
It illustrates how, for this case, the initial guess is nearly perfect if the extra term
is employed, while it is quite poor without it.
In particular, the relative residuals are as low as $10^{-5}$ - $10^{-9}$,
whereas the relative residuals for the initial guess using \cite{Zeman_2017} are in the order of
$10^{-2}$ - $10^{-5}$.

\begin{table}[tbp]
    \centering
    \caption{Material parameters as assumed for the ferrite and martensite phases.}
    \label{tbl:mat_par}
    \begin{tabular}{l|lllllccc}
        \hline
            Parameter &
            $E$ [GPa] &
            $\nu$ [-] &
            $\dot{\gamma}_0$ [1/s] &
            $m$ [-] &
            $\sigma_0$ [MPa] &
            \multicolumn{3}{c}{$h$ [MPa]}
        \\
        \hline
          &
          &
          &
          &
          &
          &
          \multicolumn{1}{l}{Hard.} &
          \multicolumn{1}{l}{Perf.pl.} &
          \multicolumn{1}{l}{Soft.}
        \\
            Ferrite &
            206.824 &
            0.3 &
            0.001 &
            0.05 &
            425 &
            940 &
            0 &
            -940
        \\
            Martensite &
            206.824 &
            0.3 &
            0.001 &
            0.05 &
            1180 &
            1740 &
            0 &
            -1740
        \\
        \hline
    \end{tabular}
\end{table}

\begin{figure}[htbp]
    \begin{subfigure}[t]{0.40\linewidth}
        \centering
        \includegraphics[width=1.\linewidth]{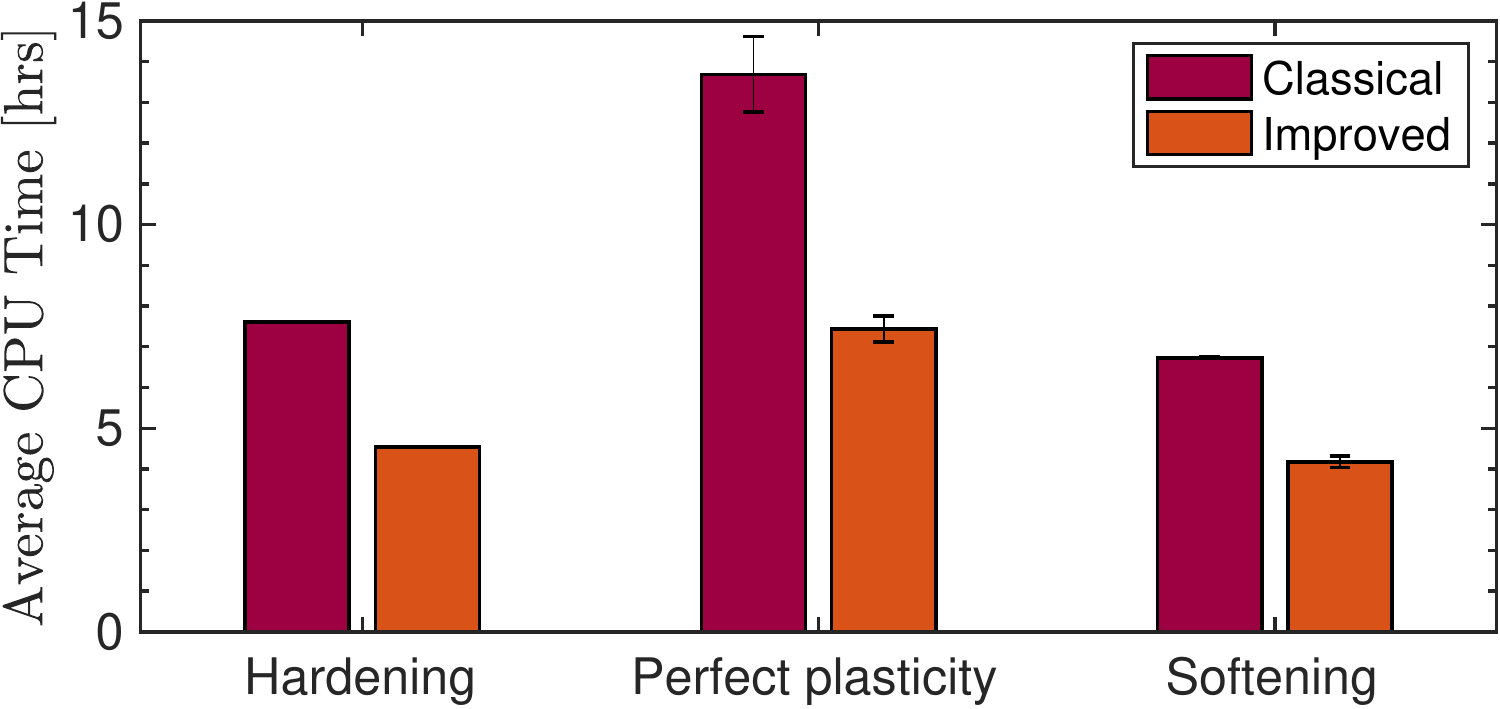}
        \caption{Total CPU time averaged over three runs for different hardening parameters.}
        \label{fig:res_scen_2_a}
    \end{subfigure}
    \hfill
    \begin{subfigure}[t]{0.55\linewidth}
        \centering
        \includegraphics[width=1.\linewidth]{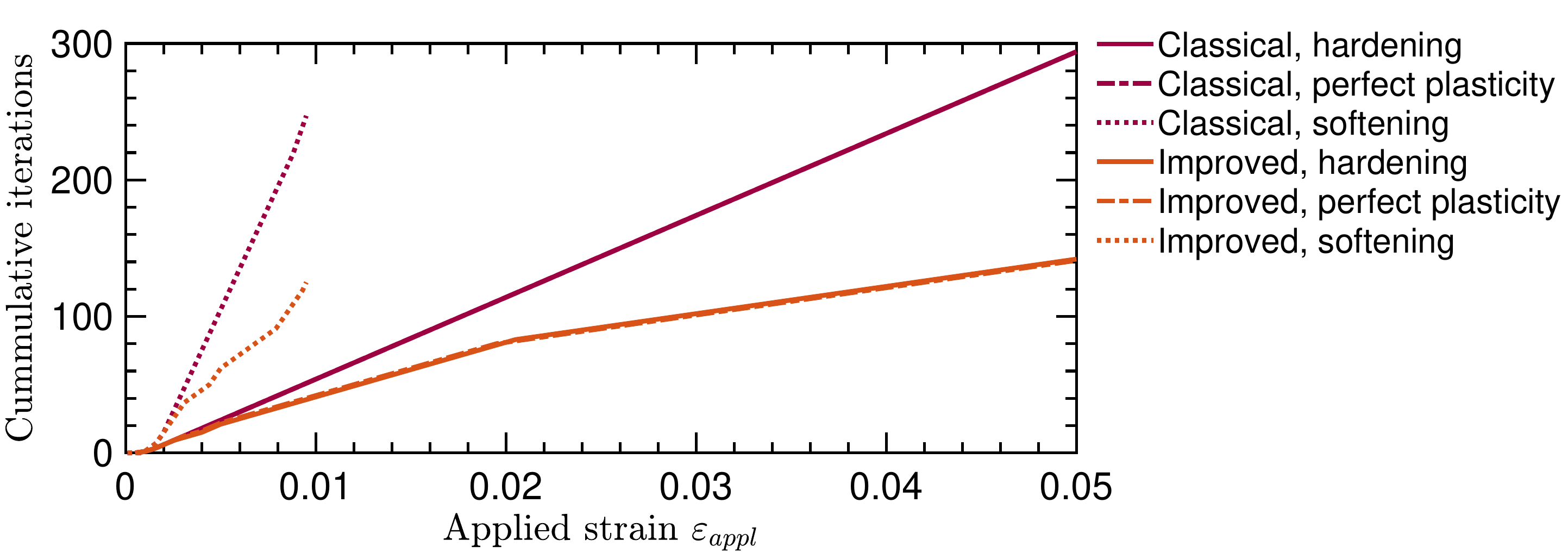}
        \caption{
            Cumulative number of Newton-Raphson iterations versus applied strain
            $\varepsilon_{appl}$ with overlapping curves for hardening and perfect plasticity.}
        \label{fig:res_scen_2_b}
    \end{subfigure}
    \caption{
        Comparison of the numerical performance of the initial guess from \cite{Zeman_2017}
        (\textit{classical}) and the initial guess from this work (\textit{improved})
        for visco-plasticity with hardening, perfect plasticity and softening.}
    \label{fig:res_scen_2}
\end{figure}

The average CPU time used to compute the mechanical response of the microstructure is shown in
Fig.~\ref{fig:res_scen_2_a}.
For each considered case, the extra term decreases the CPU time by approximately 45\%,
by reducing the number of Newton-Raphson iterations per load increment.
In particular since the initial guess is closer to the final solution,
the convergence of the Newton-Raphson protocol is improved.
This is confirmed by the convergence of the relative residual norm in Table \ref{tab:NR_conv}.
As the extra term essentially calculates the increase in plastic strain using the plastic strain
rate from the previous time step, it is most accurate where there is little change in plastic
strain rate between different time steps, thus especially for the steady state regime of the
perfectly visco-plastic model.
This is further illustrated in Fig.~\ref{fig:res_scen_2_b}, which shows the cumulative number of
iterations as a function of the imposed strain.
As observed, the improvement becomes clearly visible once the yield stress of the material is
reached and the plastic strain starts to develop at a significant rate.
For the overlapping curves of the hardening and perfectly plastic case,
the effect of the improvement is observable from an
applied strain of $\varepsilon_{appl} = 0.005$ onwards.
At $\varepsilon_{appl} = 0.02$, the predictability of the plastic strain increases and the
improvement distinguishes itself even further.
Due to the more localised nature of the plastic strain as the result of softening, the improvement
in terms of the number of Newton-Raphson iterations is not monotonic.

\begin{figure}[htbp]
    \begin{subfigure}[t]{0.45\linewidth}
        \centering
        \includegraphics[width=1.\linewidth]{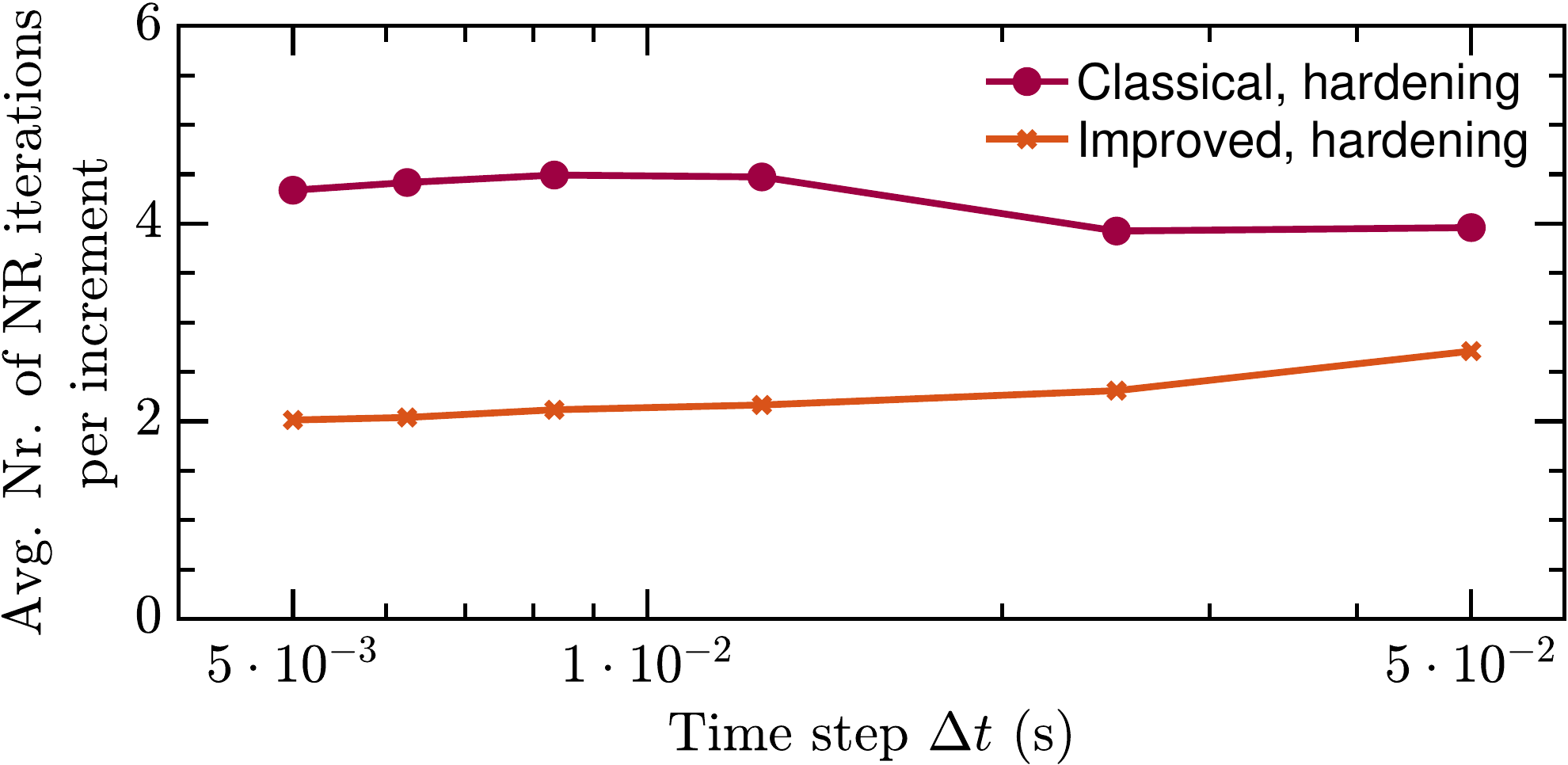}
        \caption{Influence of the time-step $\Delta t$.}
        \label{fig:avg_NR_its_dt}
    \end{subfigure}
    \hfill
    \begin{subfigure}[t]{0.45\linewidth}
        \centering
        \includegraphics[width=1.\linewidth]{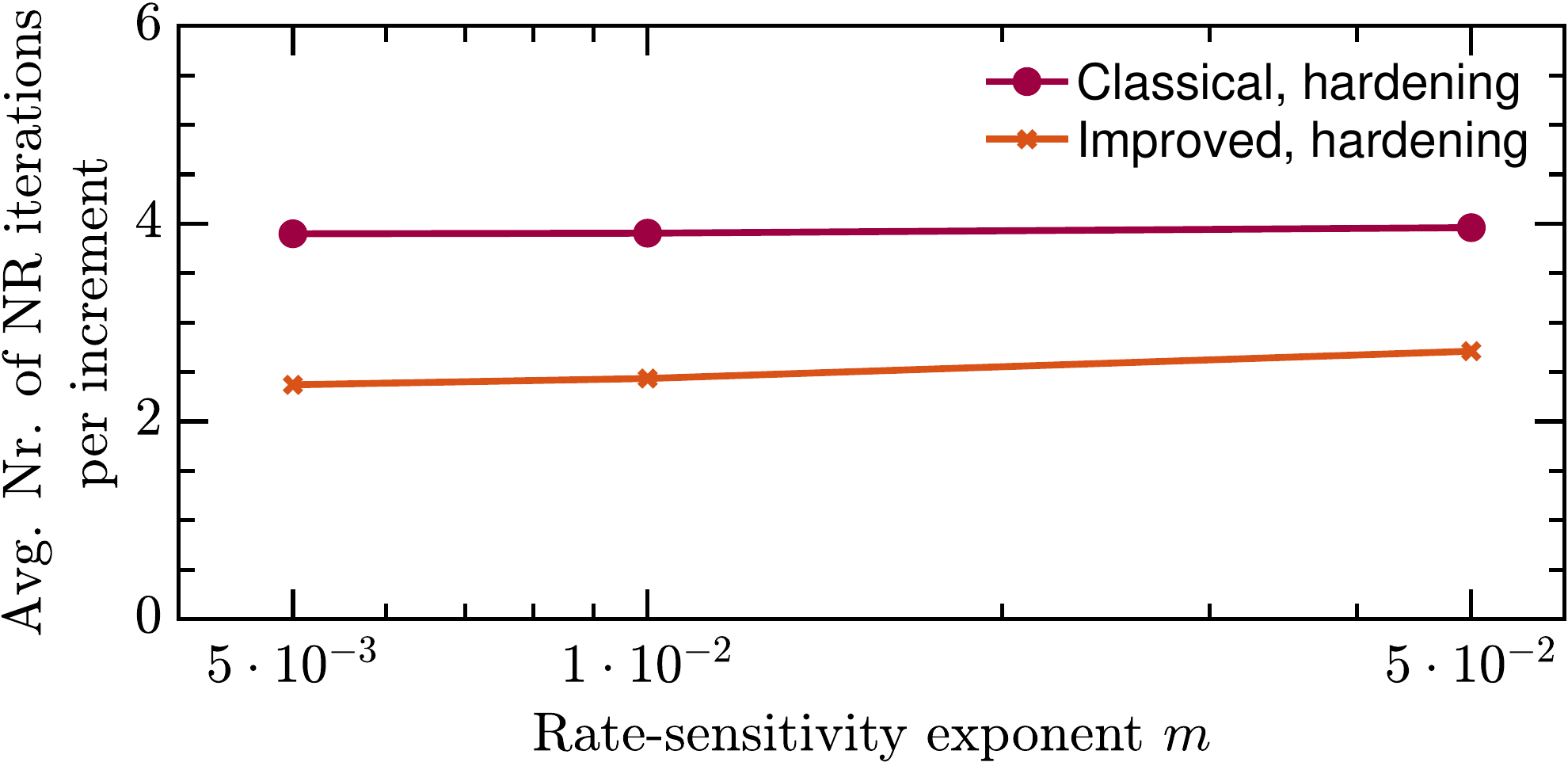}
        \caption{Influence of the rate-sensitivity exponent $m$.}
        \label{fig:avg_NR_its_m}
    \end{subfigure}
    \caption{
        Comparison of the average number of Newton-Raphson iterations per time increment of
        the initial guess from \cite{Zeman_2017} (\textit{classical})
        and the initial guess from this work (\textit{improved})
        for visco-plasticity with hardening,
        measured on a 101x101 section of the SEM image of Fig.~\ref{fig:micrograph_a}.
    }
    \label{fig:avg_NR_its}
\end{figure}

When varying the time-step size $\Delta t$ or the rate-sensitivity exponent $m$,
as shown in Fig.~\ref{fig:avg_NR_its},
the average number of Newton-Raphson iterations used per load increment is
consistently halved by the improvement.
Note that the number of time steps was increased with decreasing rate-sensitivity exponent $m$,
 such that the ratio between the two was kept constant.

\begin{table}[htbp]
    \centering
    \caption{
        Convergence of the norm of the mechanical equilibrium equation normalised with the
        yield stress of martensite during Newton-Raphson iterations for the classical and
        improved initial guess, measured at $\varepsilon_{appl} = 0.05$ for the hardening
        and perfectly plastic case, and at $\varepsilon_{appl} = 0.0075$ for softening,
        i.e.\ well in the plastic regime.}
    \label{tab:NR_conv}
    \resizebox{.6\textwidth}{!}{
        \begin{tabular}{@{}l|llllll@{}}
        \toprule
            Iteration &
            \multicolumn{2}{l}{Hardening} &
            \multicolumn{2}{l}{Perfect plasticity} &
            \multicolumn{2}{l}{Softening}
        \\
            &
            Classical &
            Improved &
            Classical &
            Improved &
            Classical &
            Improved
        \\
        \midrule
            $i = 0$ &
            4.86E+00 &
            5.86E-04 &
            4.64E+00 &
            2.77E-04 &
            3.26E+00 &
            1.42E-03
        \\
            $i = 1$ &
            1.22E-01 &
            1.49E-08 &
            1.62E-01 &
            1.35E-08 &
            1.36E-01 &
            3.30E-07
        \\
            $i = 2$ &
            4.96E-04 &
            - &
            6.87E-04 &
            - &
            2.02E-03 &
            -
        \\
            $i = 3$ &
            2.47E-08 &
            - &
            1.01E-07 &
            - &
            1.58E-06 &
            -
        \\
        \bottomrule
        \end{tabular} %
    }
\end{table}

\section{Synopsis}

A general linearisation procedure for the consistent tangent of a small-strain
visco-plastic material model was presented in this note.
The procedure is based on multi-variable linearisation around a so-called ``reference state''.
In particular, the linearisation of the time integration scheme (i.e.\ Eq.~\eqref{eq:lin_BE})
was found to yield an extra term compared to classical expressions
\cite{Ju_1990, Souza_2011, Simo_2006}, which only appears because the
material response is time-dependent.
It has the effect of yielding a very accurate initial guess for the Newton-Raphson protocol
based on the ongoing viscous flow.
It was shown, using a modern variational FFT-based solver,
that the extra term reduces both the CPU time and the number of Newton-Raphson iterations
by around a factor two.

\section*{Acknowledgement}

T.G.~was partly financially supported by
The Netherlands Organisation for Scientific Research (NWO) by a NWO Rubicon grant
number 680-50-1520.
J.V.~gratefully acknowledges Jan Zeman, Jaroslav Vond\v{r}ejc and Luv Sharma for their
input in discussions concerning FFT-based spectral methods.

\bibliographystyle{unsrtnat}
\bibliography{library}

\appendix

\section{Nomenclature}
\label{app:nomenclature}

We use boldface symbols to denote
vectors $\bm{a} = a_i \, \vec{e}_i$,
second-order tensors, $\bm{A} = A_{ij} \, \vec{e}_i \vec{e}_j$, and
fourth-order tensors, ${}^4\bm{A} = A_{ijkl} \, \vec{e}_i \vec{e}_j \vec{e}_k \vec{e}_l$.
A tensor contraction is denoted using centered dot,
e.g.\ $\bm{C} = \bm{A} \cdot \bm{B}$ corresponds to $C_{ik} = A_{ij} B_{jk}$.
A double tensor contraction is denoted using colon, e.g.\ $c = \bm{A} : \bm{B}$ corresponds to
$c = A_{ij} B_{ji}$.
$\bm{I} \equiv \delta_{ij} \, \vec{e}_i \vec{e}_j$ is a second-order unit tensor,
and
$\bm{I} \bm{I} \equiv \delta_{ij} \delta_{kl} \, \vec{e}_i \vec{e}_j \vec{e}_k \vec{e}_l$
corresponds to a dyadic product of two second-order unit tensors.
$\text{tr}(\bm{A}) \equiv A_{ii} / 3$ is the trace of second-order tensor.
$\nabla \cdot \bm{A}$ corresponds to the divergence operator $\partial A_{ij} / \partial x_j$.
Note that for all index notations a summation of the three spatial dimensions is implied.

\section{Implementation of the improved initial guess}
\label{app:implementation}

We use Algorithm 1 of Ref.~\cite{Zeman_2017} whereby the only modification is line 7
of the algorithm, which now reads
\begin{equation}
    \label{eq:improved_ig}
    \underline{{{\bm{{G}}}}} \, \underline{{{\bm{C}}}}_{\left(t\right)} \,
    \delta {\underline{\bm{\varepsilon}}}^{*}_{\left(0\right)}
    =
    -\underline{{{\bm{{G}}}}} \, \underline{{{\bm{C}}}}_{\left(t\right)}
    \,\left[{ \bm{{\underline{E}}}}_{\left(t+\Delta t \right)} -
    \bm{{\underline{E}}}_{\left(t\right)} -
    {\Delta t{{\dot {\underline{\gamma }}^t}}}\underline{\kappa}^t
    \,{\bm{\underline{N}}^t} \right]
\end{equation}
see Ref.~\cite{Zeman_2017} for nomenclature. Here we only specify that we take all nodes
(grid points) visco-plastic and that
$\dot{\underline{\gamma }}^t$, $\underline{\kappa}^t$ and ${\bm{\underline{N}}^t}$
are columns that collect the nodal quantities.
Note also that we use Algorithm 1 of Ref.~\cite{Zeman_2017},
without any modification, as reference.

\section{Generalised trapezoidal integration}
\label{app:theta_integration}

We now generalise our results to the generalised trapezoidal integration scheme,
which employs a linear combination of variables evaluated at time $t$ and at time $t + \Delta t$
through a parameter $0 \leq \theta \leq 1$.
Note that the choice of this parameter allows one to recover the explicit
forward Euler scheme when $\theta = 0$ and the backward Euler scheme when $\theta = 1$.
The drawback of the generalised trapezoidal scheme lies in its return map,
which requires that the following set of non-linear equations is solved
\begin{equation}
    \begin{split}
        \bm{\varepsilon}^{t+\Delta t}_e &=
        {}^{tr}\bm{\varepsilon}_e - \Delta \gamma
        \left[ \left( 1 - \theta \right) \bm{N}^t + \theta \bm{N}^{t+\Delta t} \right]
        \\
        \Delta \gamma &=
        \Delta t \left(
            \left( 1 - \theta \right) \dot{\gamma}^t +
            \theta \dot{\gamma}^{t+\Delta t}
        \right)
    \end{split}
\end{equation}
In comparison, the implicit backward Euler scheme ($\theta = 1$)
only requires the solution of the latter, non-linear scalar, equation for $\Delta \gamma$.
Using the generalised trapezoidal integration scheme,
Eq.~\eqref{eq:epspdot} is discretised as follows
\begin{equation}
    \label{eq:Delta_epsp_trapz}
    \Delta \bm{\varepsilon}_{p} =
    \Delta \gamma \; \underbrace{\left(\left( 1 - \theta \right) \bm{N}^t +
    \theta \bm{N}^{t+\Delta t} \right)}_{\bm{N}^{\theta}}
\end{equation}
where
\begin{equation}
    \label{eq:dgam_trapz}
    \Delta \gamma = \gamma^{t+\Delta t} - \gamma^t
    = \left( 1 - \theta \right) \Delta t{{\dot \gamma}^t}
    + \theta \Delta t{{\dot \gamma }^{t + \Delta t}}
\end{equation}
After linearising Eq.~\eqref{eq:Delta_epsp_trapz} around ${\bullet^*}$, we obtain
\begin{equation}
    \label{eq:del_epsp_trapz}
    \delta \bm{\varepsilon}_p =
    \delta \gamma \, \bm{N}^\theta + \underbrace{
        \left( \gamma^* - \gamma^{* - \Delta t} \right)
    }_{\Delta \gamma^*} \theta \delta \bm{N}
\end{equation}
where the quantities $\bm{N}^\theta$ and $\Delta \gamma^*$ are fully known.
The expression for $\delta \bm{N}$ is derived by linearising its definition
in Eq.~\eqref{eq:N} around ${\bullet^*}$ as
\begin{equation}
    \label{eq:del_N_trapz}
    \delta \bm{N} =
    \left[ \frac{3G}{\sigma_{eq}^*} {}^4 \bm{I}^d
    - \frac{2G}{\sigma_{eq}^*} \bm{N}^* \bm{N}^* \right]
    : \delta \bm{\varepsilon}_e
\end{equation}
To find $\delta \gamma$, we combine Eqs.~\eqref{eq:dgam_trapz}, \eqref{eq:gamma_lin} and
\eqref{eq:gamma_dot_lin} into
\begin{equation}
    \label{eq:lin_dgam}
    \gamma^* - \gamma^t - \left( 1 - \theta \right) \Delta t \dot{\gamma}^t -
    \theta \Delta t \dot{\gamma}^*
    + \delta \gamma - \theta \Delta t \delta\dot \gamma = 0
\end{equation}
The small variation $\delta\dot \gamma$ is then derived as
\begin{equation}
    \label{eq:lin_del_gam_trapz}
    {\delta\dot \gamma } =
    \frac{\partial \dot{\gamma}}{\partial \sigma_{eq}} \delta \sigma_{eq}
    + \frac{\partial \dot{\gamma}}{\partial \sigma_s} \delta \sigma_s =
    \frac{\alpha^*}{\Delta t} \left( 2 \bm{N}^* : \delta \bm{\varepsilon}_e -
    \frac{\sigma_{eq}^* \, h \, \delta \gamma}{\sigma_s^* G} \right)
    ,
    \qquad
    \alpha^* =
    \frac{\dot{\gamma}_0 G \Delta t}{m \sigma_s^*}
    \left( \frac{\sigma_{eq}^*}{\sigma_s^*} \right)^{\frac{1}{m} - 1}
\end{equation}
in which $\alpha^*$ is unaffected by the choice of integration scheme.
A closed form expression for $\delta \bm{\varepsilon}_p$ can now be obtained from
Eqs.~(\ref{eq:dgam_trapz}--\ref{eq:lin_dgam}).
This step of the procedure exposes the major disadvantage of
the generalised trapezoidal integration scheme.
The lack of co-linearity between ${}^{tr}\bm{N}$ and $\bm{N}^*$ requires a system
of non-linear equations to be solved during the return map,
as opposed to the single non-linear scalar equation for the backward Euler scheme.
As a result, both the linearised stress update:
\begin{equation}
    \label{eq:del_sigma_a_trapz}
    \delta \bm{\sigma} =
    {}^4{\bm{C}_{vp}^*}: \delta \bm{\varepsilon} -
    \frac{G \beta^*}{\alpha^*}{\left( {{\theta \Delta t{{\dot \gamma }^{*}}
    + \left( 1 - \theta \right) \Delta t{{\dot \gamma }^{t}} - \gamma^* + \gamma^t}} \right)}
    \, \, \left({}^4{\bm{P}^{*}}\right)^{-1}:{\bm{N}^\theta}
\end{equation}
and the consistent tangent:
\begin{equation}
    \label{eq:gen_C_vp_trapz}
    {}^4{\bm{C}_{vp}^*} =
    {}^4{\bm{C}_e} -
    2G \theta \beta^* \left({}^4{\bm{P}^{*}}\right)^{-1}:{\bm{N}^{\theta}}{\bm{N}^*} -
    \frac{\Delta {\gamma^*} 4G^2 \theta}{\sigma_{eq}^*} \left({}^4{\bm{P}^{*}}\right)^{-1} :
    \left[ \tfrac{3}{2} {}^4 \bm{I}^d - \bm{N}^* \bm{N}^* \right]
\end{equation}
contain the inverse of the fourth-order tensor
\begin{equation}
    \label{eq:P_trapz}
    {}^4 \bm{P}^{*} =
    \left( 1 + \frac{3 G \theta \Delta \gamma^*}{\sigma_{eq}^*} \right) {}^4 \bm{I}^d -
    \frac{2 G \theta \Delta \gamma^*}{\sigma_{eq}^*} {\bm{N}^*}{\bm{N}^*}
\end{equation}
The constant $\beta^*$ reads
\begin{equation}
    \label{eq:beta_trapz}
    \beta^* = \frac{2 \alpha^*}{1 + 2 \alpha^* \theta \, \bm{N}^{\theta} : \bm{N}^* +
    \frac{\sigma_{eq}^* h}{\sigma_s^* G} \theta \alpha^*}
\end{equation}
Note that in the equation above the implicit backward Euler scheme is recovered for $\theta = 1$
as the product~$\bm{N}^{\theta}:{\bm{N}^*} = \frac{3}{2}$ after applying the linearisation as
defined in Eq.~\eqref{eq:N_lin}.

\begin{figure}[htpb]
    \centering
    \includegraphics[width=.45\linewidth]{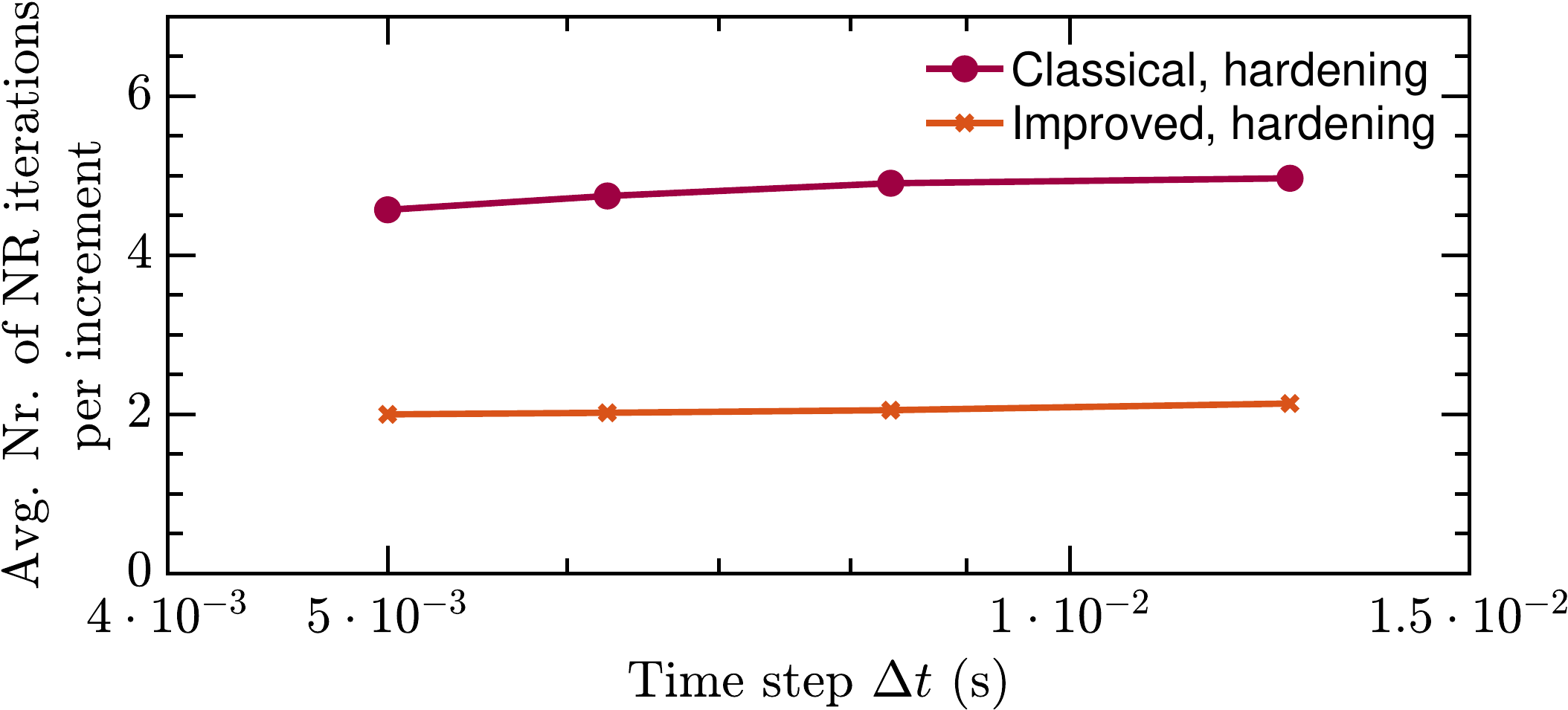}
    \caption{
        Comparison of the average number of Newton-Raphson iterations versus the time-step
        $\Delta t$ of the initial guess from \cite{Zeman_2017} (\textit{classical})
        and the initial guess from this work (\textit{improved}),
        using visco-plasticity with hardening and the generalised
        trapezoidal integration scheme with $\theta = 0.5$.
        These measurements were performed on a 101x101 section of the
        SEM image of Fig.~\ref{fig:micrograph_a}.
    }
    \label{fig:NR_its_trapz}
\end{figure}

Similar to the main text, we have two relevant choices for the reference state $\bullet^*$.
An `ordinary' Newton-Raphson iteration is recovered by taking
${\bullet^*} \equiv {\bullet^{i}}$ at $t + \Delta t$.
In this case Eq.~\eqref{eq:del_sigma_a_trapz} reduces to the classical
\begin{equation}
    \label{eq:del_sigma_NR_trapz}
    \delta \bm{\sigma} = {}^4 \bm{C}_{vp}^{i} :  \delta \bm{\varepsilon}
\end{equation}
This result follows from combining the discretised strain rate Eq.~\eqref{eq:dgam_trapz} with
Eq.~\eqref{eq:del_sigma_a_trapz}, cancelling all four terms within brackets in the latter.
The improved initial guess is found by taking ${\bullet^*} \equiv {\bullet^{t}}$.
In this case $\gamma^* \equiv \gamma^t$ reducing the terms within brackets of
Eq.~\eqref{eq:del_sigma_a_trapz} to $\Delta t \dot \gamma^t$.
After some reorganisation, we can write the final result for the stress update in a form
identical to the main text, namely
\begin{equation}
    \label{eq:del_sigma_c_trapz}
    \delta \bm{\sigma} =
    {}^4{\bm{C}_{vp}^t}:\left[ {\delta \bm{\varepsilon} -
    \Delta t \dot{\gamma}^t \kappa^t {\bm{N}^t}} \right]
\end{equation}
in which the change of integration scheme is only observed by a small change in the constant
\begin{equation}
    \kappa^t =
    \frac{1}{\left( 1 + \frac{{\sigma_{eq}^*}h}{{\sigma_{s}^*} G} \theta \alpha^* \right)}
\end{equation}
Indeed, for $\theta = 1$ Eq.~\eqref{eq:kappa} is recovered.
Fig.~\ref{fig:NR_its_trapz} shows the average number of Newton-Raphson iterations used
per load increment for different time-step sizes $\Delta t$ using the trapezoidal scheme
with $\theta = 0.5$.
The improved initial guess gives results consistent with that of the backward Euler schemes as
the Newton-Raphson iterations are approximately halved, saving significant CPU-time.

\end{document}